\documentclass[12pt,letterpaper]{article}
\usepackage{subcaption}
\usepackage[a4paper, total={7in, 10in}]{geometry}

\usepackage{graphicx}
\usepackage{diagbox}
\usepackage{helvet}
\usepackage{authblk}
\usepackage{hyperref}
\usepackage{amsmath} 
\usepackage{amssymb} 
\usepackage{orcidlink} 
\usepackage[super,comma,sort&compress]  
   {natbib}
\usepackage[right]{lineno} \linenumbers
\usepackage{algorithm}
\usepackage{algpseudocode}
\usepackage{textcomp}
\usepackage{xcolor}
\usepackage{longtable}
\usepackage{hyperref}
\usepackage{multirow}
\usepackage{pifont}
\usepackage{booktabs}
\usepackage{enumitem}
\usepackage[dvipsnames]{xcolor}
\usepackage[normalem]{ulem}

\newcommand{\hdk}[1]{\textcolor{black}{#1}}

\makeatletter
\renewcommand{\maketitle}{\bgroup\setlength{\parindent}{0pt}
\begin{flushleft}
  \textbf{\@title}
  
  \@author
\end{flushleft}\egroup}
\makeatother

\title{\hdk{Multi-agent Power Grid Restoration Under Uncertainty Considering Coupled Transportation-Power Networks}}
\date{}

\author[1,\orcidlink{0000-0002-1159-5084}]{Harshal D. Kaushik}
\author[1,2]{Roshni Anna Jacob}
\author[3]{Souma Chowdhury}
\author[1,2,*,\orcidlink{0000-0001-7478-5670}]{Jie Zhang}

\affil[1]{Department of Mechanical Engineering, The University of Texas at Dallas, Richardson, TX 75080, USA.}
\affil[2]{Department of Electrical and Computer Engineering, The University of Texas at Dallas, Richardson, TX 75080, USA.}

 \affil[3]{Department of Mechanical and Aerospace Engineering, University at Buffalo, Buffalo, NY 14260, USA.}

\affil[*]{Correspondence: jiezhang@utdallas.edu}

\begin{document}

\maketitle

\section*{SUMMARY}

\hdk{Restoring power distribution systems after extreme events such as tornadoes presents significant logistical and computational challenges. The complexity arises from the need to coordinate multiple repair crews under uncertainty, manage interdependent infrastructure failures, and respect strict sequencing and routing constraints. Existing methods often rely on deterministic heuristics or simplified models that fail to capture the interdependencies between power and transportation networks, do not adequately model uncertainty, and lack representation of the interrelated dynamics and dependencies among different types of repair crews--leading to suboptimal restoration outcomes. To address these limitations, we develop a stochastic two-stage mixed-integer programming framework for proactive crew allocation, assignment, and routing in power grid restoration. The primary objective of our framework is to minimize service downtime and enhance power restoration by efficiently coordinating repair operations under uncertainty. Multiple repair crews are modeled as distinct agents, enabling decentralized coordination and efficient task allocation across the network. To validate our approach, we conduct a case study using the IEEE 8500-node test feeder integrated with a real transportation network from the Dallas-Fort Worth (DFW) region. Additionally, we use tornado event data from the DFW area to construct realistic failure scenarios involving damaged grid components and transportation links. Results from our case study demonstrate that the proposed method enables more coordinated and efficient restoration strategies. The model facilitates real-time disaster response by supporting timely and practical power grid restoration, with a strong emphasis on interoperability and crew schedule coordination.} 

\section*{KEYWORDS}


Power restoration, Transportation network, Uncertainty, Two-stage stochastic optimization, Resource allocation and assignment, Optimal crew deployment.

\section{INTRODUCTION}

\hdk{Electricity is the lifeblood of modern society. When power is disrupted, so is the stability of everything we depend on, including transportation, communication, healthcare, and emergency services. Extreme weather events are among the primary causes of such disruptions, often leading to cascading failures across critical infrastructure and leaving communities immobilized for days or even months. The recovery process becomes even more challenging when physical damage impacts both the power grid and its interconnected infrastructure systems. Restoring power distribution networks under these conditions presents a significant computational challenge. These networks extend across vast regions, and the type and extent of damage are inherently uncertain. This uncertainty, combined with the grid’s reliance on transportation networks for repair crew deployment, makes restoration planning a highly complex and interconnected problem. To address this large-scale unpredictable problem, and reduce the network downtime, we need approaches to abstract and/or decompose the problem into forms that can be solved efficiently. In this work, we seek to address this by solving resource planning ahead of time over different scenarios, generated by integrating transportation, power grid, and extreme event data. Additionally, we account for the logistical challenge of deploying multiple crews across the network.}

Extreme weather events have become the primary cause of power outages in the U.S. power grid \cite{Parera_climate, Blackout_Kenward}. These outages have significantly disrupted societal operations and resulted in substantial economic losses \cite{Resilience_assessment_Tari}. For instance, the February 2021 Texas winter storm, which left millions without power, incurred estimated costs ranging from 80 billion to 130  billion. In 2017, Hurricane Harvey and Hurricane Irma resulted in power outages affecting nearly 300,000 \cite{Harvey_hampers} and 15 million customers \cite{Imma_hampers}, respectively. Consequently, enhancing outage management and expediting network restoration are critical tasks for utility companies after a disastrous event.

To enhance the power grid resilience in response to extreme events, various measures have been proposed in the literature, which can be broadly classified into mainly three categories: pre-event, during-event, and post-event measures \cite{Shuai_post_storm}. Pre-event strategies, such as outage prediction and damage estimation~\cite{Pred_Hurricane}, grid hardening and resource allocations such as distributed energy resources (DERs) and blackstart generators~\cite{Robust_Opt_Yuan,Black_start_resource_allo}, as well as proactive reconfiguration and generation dispatch~\cite{Amirioun_resilience_practive_plann}, help system operators mitigate impacts and costs by preparing in advance. During-events strategies include outage detection~\cite{ChenRoshni2023PS}, network reconfiguration and emergency load shedding~\cite{ZhangRoshniSteveChowdhury2023}, intentional islanding~\cite{SobhanRoshniIntentionalIsland2023}, etc. Post events, the utilities assess the actual damages and restore the power using system repair, generator start up, and load restoration strategies \cite{Milestones_Hou, Optimal_Generator, linear_topological_constraints_Shi, SOCP}.

Restoration in power networks include service and infrastructure restoration~\cite{force2022methods}. Service restoration includes actions such as network reconfiguration and intentional islanding to mitigate failures and support critical loads~\cite{ZhangRoshniSteveChowdhury2023}. While this form of service restoration improves the response to emergency events, it however does not bring the network back to its pre-disrupted state. This can only be resolved by repairing the damaged components, which involves resource planning, considering logistical challenges and transportation constraints. In the literature, the repair and restoration problem in power systems have been addressed by exploring different strategies such as pre-hurricane crew staging and post hurricane resource management~\cite{ArabKhodaeiHanKhator2015}, multi-stage approaches for crew routing considering network power-flow \cite{HentenryckCoffrinBent2011}, and randomized adaptive vehicle decomposition \cite{SimonCoffrinHentenryck2012}. Other methods include queuing theory for repair schedules \cite{ZapataSilvaGonzalezBurbanoHernandez2008}, dynamic programming for crew dispatch \cite{CarvalhoCarvalhoFerreira2007}, and large-scale complicated constrained optimization algorithms for crew assignments \cite{XuGuikemaDavidsonNozick2007, HDK2019, HDK2020, HDK_SIAM, HDK_TAC}. 

\hdk{Several critical factors are yet to be explored in prior works on power grid restoration. First, when scheduling repairs, it is practical to consider multiple crew requirements depending on the type and the extent of the damage on the network. A significant challenge in addressing this is, however, the inherent stochasticity of such events. Second, to the best of our knowledge, existing research has largely overlooked the integration of real-world road networks into restoration planning, despite their crucial role in routing repair crews. Moreover, extreme events can also disrupt the transportation network itself, further complicating the routing process. Accurately estimating repair times and damage extents, along with potential transportation network failures, remains highly challenging. Here, we build upon the methodology proposed in~\cite{ArifWangWangChen2018, ArifWangChenWang1, KPEC_HDK}. Our work advances the research on restoration and repair processes by considering coupled power and transportation networks under uncertainty. The proposed framework is designed to enable practical and timely restoration of the power grid, with an emphasis on the interoperability and coordination of repair schedules. The key contributions of this study are as follows:}

\begin{itemize}[leftmargin=0.36cm]
    \item \hdk{A realistic transportation network is considered alongside the power grid to reflect real-world constraints on crew mobility. Specifically, we model the primary level of the power distribution network and superimpose it on actual road network data. This integration results in a large-scale structure due to the spatial and operational complexity of coordinating repairs across coupled networks. To manage computational tractability, we simplify the problem by leveraging graph-theoretic concepts, enabling efficient representation of the network structure.}
    \item \hdk{To ensure uncertainty awareness and mitigation, we account for the inherent stochasticity involved in distribution grid restoration following a storm. Various scenarios are generated in advance by incorporating uncertainties and worst-case outcomes in repair times, repair demands, and transportation network disruptions. To ground our analysis in a realistic context, we utilize tornado impact data from the Dallas-Fort Worth (DFW) area. Using this dataset, we identify critical components of the power and transportation networks that are most vulnerable to storm-induced failures. Based on this analysis, we model realistic damage scenarios, including both damaged grid nodes and disrupted road segments.}
    \item \hdk{A two-stage stochastic mixed-integer programming model is developed to address the distribution system repair and restoration problem under uncertainty. In the first stage, the model determines the optimal allocation of repair crew resources before the storm using recourse constraints across all scenarios. Once the crew capacities and their sequence of deployment are finalized, the second stage addresses the multi-agent vehicle routing problem by optimally routing the crews to damaged nodes under each scenario. This structured approach ensures that both strategic resource planning and operational routing decisions are made with full awareness of uncertainty.}
    \item \hdk{To provide a more practical and realistic restoration schedule, we develop a detailed plan that allocates multiple specialized crews based on the severity and type of failure at each damaged node. Each crew is assigned a specific function such as line repair, transformer replacement, or vegetation clearing, and operates within a defined time window. This coordination is visualized through Gantt charts, which represent the repair timeline and emphasize the sequencing and dependencies between tasks. By explicitly modeling multiple crews and their interactions, our approach minimizes downtime, maximizes power restoration, and ensures efficient use of available repair resources.}
\end{itemize}

\hdk{The rest of the paper is organized as follows: Section \ref{sec:results} presents the scenario generation approach, describes the multi-crew framework and the coupled transportation-power network, and discusses simulation results through case studies. Section \ref{sec:discussion} includes a discussion of the proposed methodology, underlying assumptions, limitations, and potential directions for future work. Section \ref{sec:methods} details the mathematical formulation and outlines the two-stage stochastic programming framework. Finally, Section \ref{sec:conclude} concludes the paper with key takeaways.}

\section{RESULTS}\label{sec:results}

\begin{figure*}[thb!]
     \centering
    \includegraphics[width=0.9\textwidth]{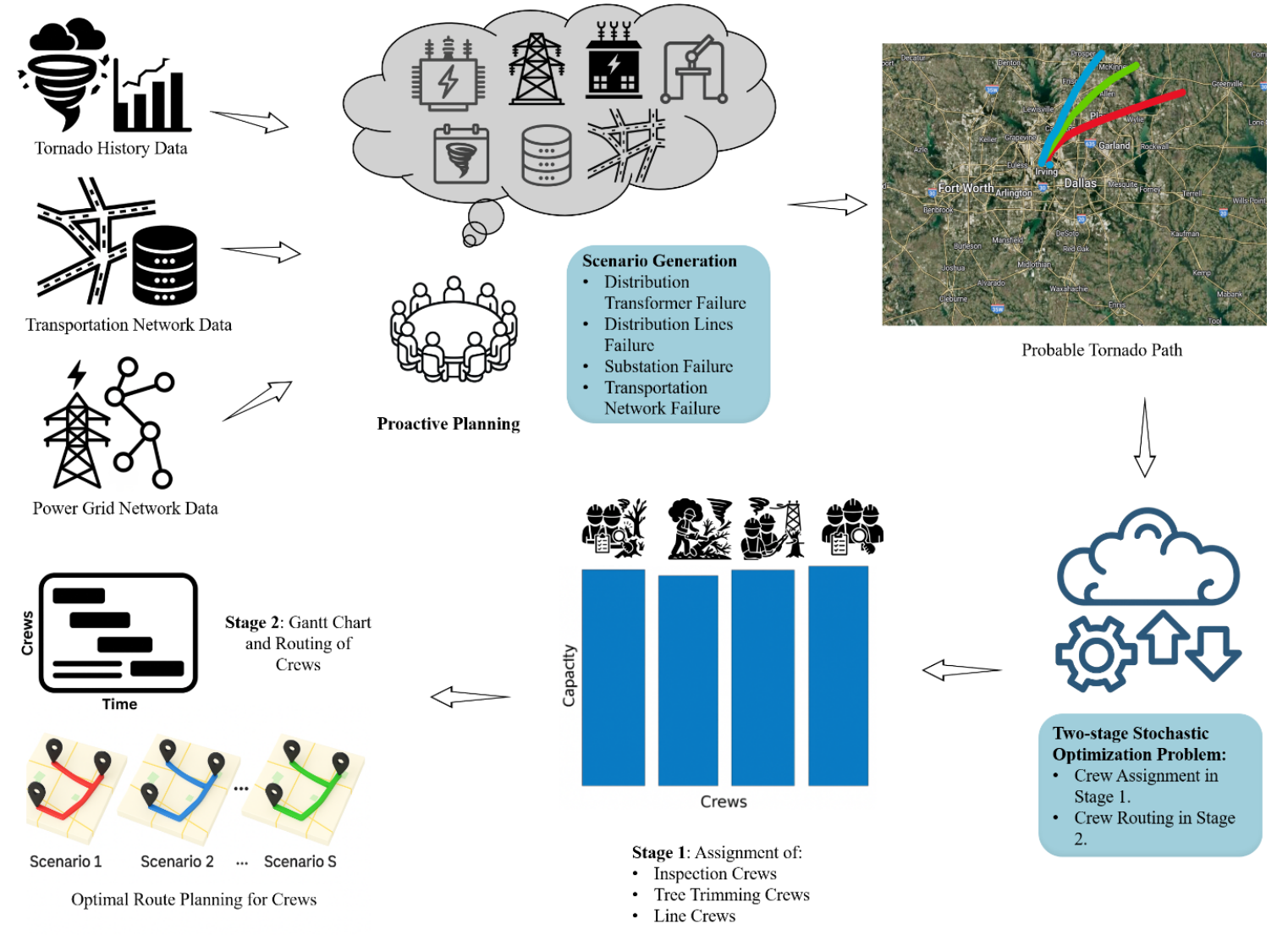}
        \caption{\hdk{Overview of the proposed process for restoring the distribution grid under uncertainty. The framework begins with scenario generation using tornado history, power grid, and transportation network data. In Stage 1, repair crew types and capacities are allocated in anticipation of damage. In Stage 2, crews are routed to damaged nodes and scheduled using Gantt charts, reflecting a realistic and coordinated multi-crew restoration plan. The approach integrates infrastructure interdependencies and stochastic failure modeling to support proactive and efficient decision-making.}}
    \label{fig:overview_restoration_process}
\end{figure*}

\subsection{Overview of Methodology}

\hdk{The primary objective of this work is to develop an uncertainty-aware scheduling framework for the efficient repair and restoration of power distribution grids, particularly under the impact of natural disasters such as tornadoes. As illustrated in Figure \ref{fig:overview_restoration_process}, the framework begins with a scenario generation process that integrates historical tornado data, transportation network topology, and power grid infrastructure to simulate a wide range of failure scenarios, including transformer outages, distribution line failures, and road blockages. These scenarios are then used as inputs to a two-stage stochastic mixed-integer optimization model. In the first stage, optimal crew capacities and assignments are determined prior to the storm, considering projected damage and resource constraints. In the second stage, once actual damage is observed, specialized crews, such as inspection, line repair, and tree trimming crews are routed and scheduled accordingly. The overall objective of the model is to assign and route repair crews in a way that maximizes power restoration while minimizing the repair time, travel time, and coordination inefficiencies. Figure \ref{fig:overview_restoration_process} highlights this process flow, showing data fusion, crew capacity planning, and post-storm routing, with the output visualized through Gantt charts and optimized scenario-specific routing plans.}



\subsection{Scenario Generation} \label{subsec:scenario_gen}
\begin{figure}[htpb!]
     \centering
    \includegraphics[width=0.81\textwidth]{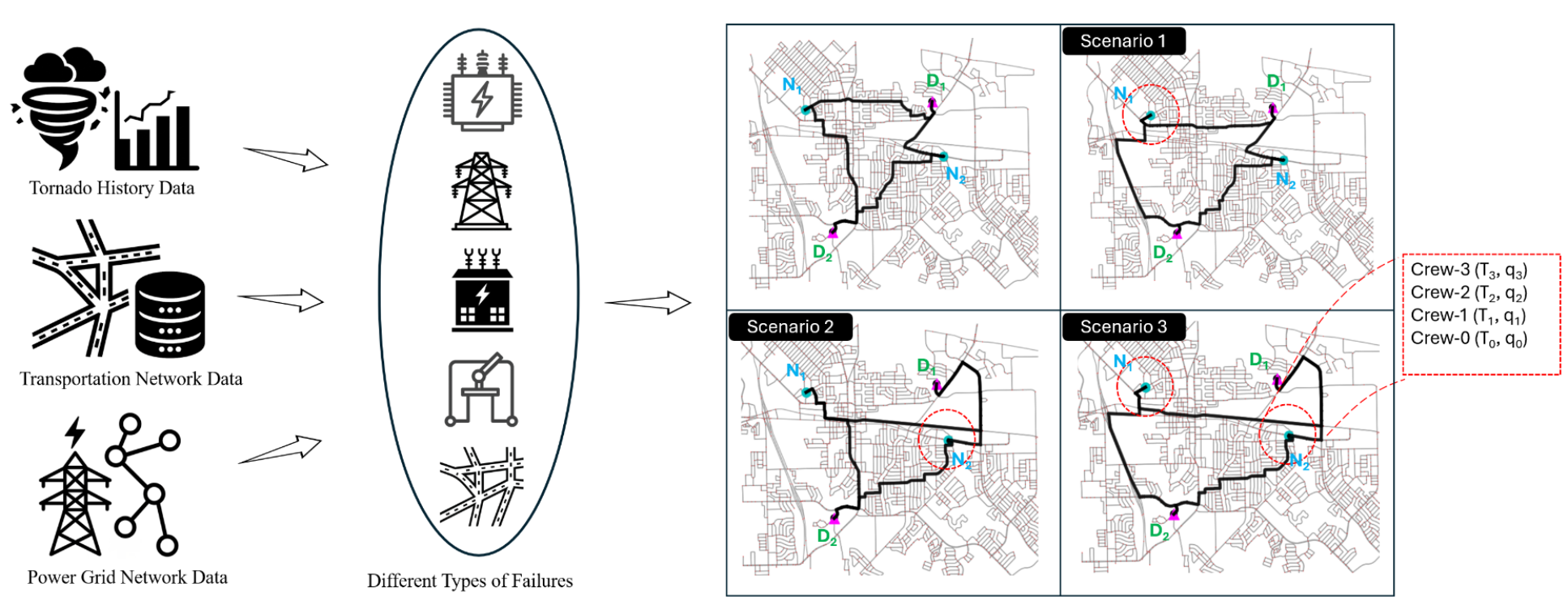}
    \caption{\hdk{Scenario generation based on variations in repair times ($\mathcal{T}$), repair demands ($\mathcal{D}$), and transportation network availability. Scenarios are derived using historical tornado data, power grid infrastructure (including transformers, substations, switches, and distribution lines), and transportation network data. Each scenario captures different combinations of failures to reflect realistic post-disaster conditions and their impact on optimal crew routing from depots ($\mathbb{D}$) to damaged nodes ($\mathbb{N}$).}}
    \label{fig:scenarios}
\end{figure}

Scenarios are generated proactively before a disaster using a rapid and coordinated approach to ensure efficient resource allocation and preparedness for a range of outcomes, including worst-case events. This process begins with damage estimation and forecasting that incorporates actual transportation network topology, historical tornado event data, and the structure of the electrical grid. Continuous monitoring and historical records of similar disasters help identify vulnerable areas, particularly those with critical power infrastructure. \hdk{The generated scenarios capture uncertainties related to transportation network failures, varying levels of damage across the power grid, and corresponding repair times. Specifically, we consider potential damage to transformers, power lines, substations, and switches, ensuring a comprehensive representation of grid vulnerabilities.}

\hdk{In our current framework, we assume these scenarios are created a priori, and we then focus on solving the repair and restoration problem based on these inputs.} The scenarios are generated using probabilistic distributions: repair times ($\mathcal T$) is calculated using a lognormal distribution $\mu = -0.3072$ and $\sigma = 1.8404$ \cite{ZhuZhouYanChen}, repair demands ($\mathcal{D}$) are sampled from a uniform distribution, and transportation failures are informed by historical tornado data. Each scenario represents a unique combination of damaged network nodes and altered transportation routes, as shown in Figure \ref{fig:scenarios}. Specifically, (a) transportation route availability is varied, affecting shortest paths between depots ($\mathbb{D}$) and nodes ($\mathbb{N}$); (b) repair times reflect the estimated severity of damage at each node, while repair demand represents the required restoration effort, quantified in terms of the number of crew members needed.

\subsection{Multiple Crews}\label{subsec:multiple_crew}
\begin{figure}[htpb!]
     \centering
    \includegraphics[width=0.99\textwidth]{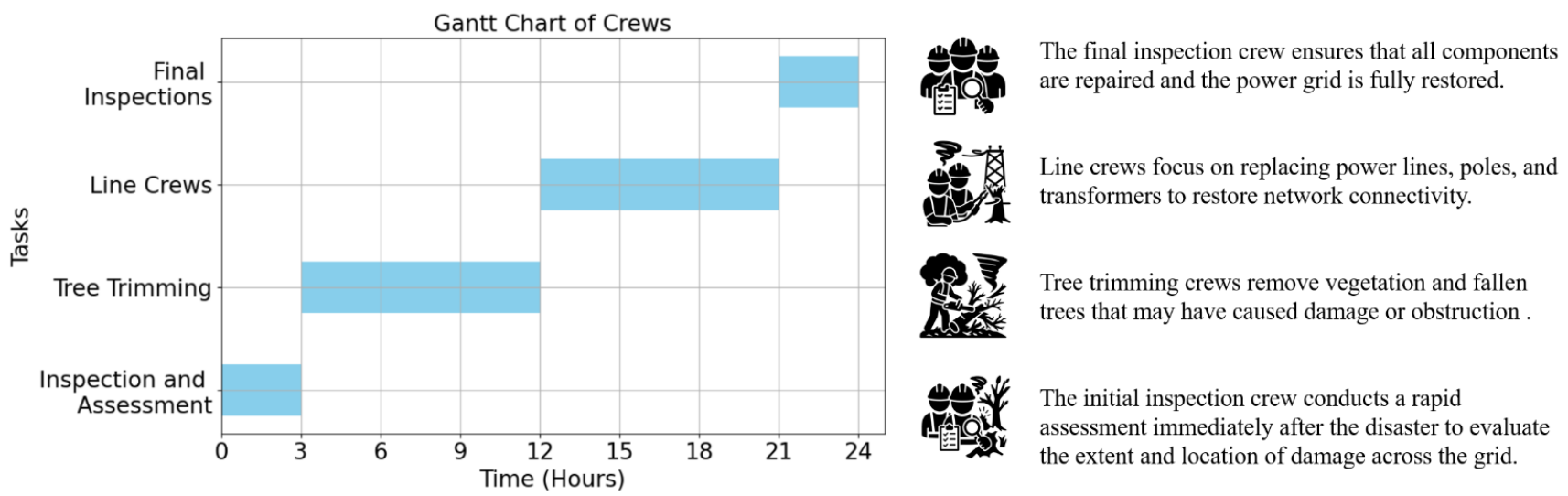}
    \caption{\hdk{Gantt chart illustrating the sequential deployment of specialized repair crews during power grid restoration. The timeline shows how inspection and assessment, tree trimming, line repairs, and final inspections are coordinated over a 24-hour period to ensure efficient and timely recovery.}}
    \label{fig:Gantt_generic}
\end{figure}

To provide a deployable repair schedule, we develop a detailed plan that assigns different crews based on the extent of failure at each damaged node. Figure \ref{fig:Gantt_generic} illustrates a sequential operation of different repair crews, responding to a particular damaged node,  where the x-axis is a tentative time taken by each crew in hours. Also,for simplicity, we do not assume overlapping between the two consecutive crews. We also consider that each crew has a specific function and operates within a defined time frame. \hdk{We consider four types of crews in this restoration process: initial inspection crews, tree trimming crews, line crews, and final inspection crews.} The inspection and assessment crew does the initial evaluation immediately after the disaster. Following this, the tree trimming crew removes any vegetation and fallen trees that may have caused or worsened the damage to the power lines. After the area is cleared, the line crew focuses on repairing and replacing damaged power lines, poles, and transformers. They perform high-voltage line work, pole climbing, and equipment operations. Finally, the process concludes with the final inspection crew. The times along x-axis in Figure \ref{fig:Gantt_generic}  are for representation. The actual time spent by each crew depends on the scale of damage to the network components.


\subsection{Integrated Transportation and Power Networks}\label{subsec:integrated_network}

To make more realistic and practical repair schedules, we solve our stochastic repair and restoration problem over the primary distribution network with critical components which is coupled with a real transportation network. This coupled transportation network can be represented as graphs consisting of nodes (depots, damaged components, intersection nodes) and edges connecting these nodes using the available road network. 
 

\begin{figure}[htpb!]
    \centering
    \begin{subfigure}[b]{0.27\textwidth}
        \centering
        \includegraphics[width=\textwidth]{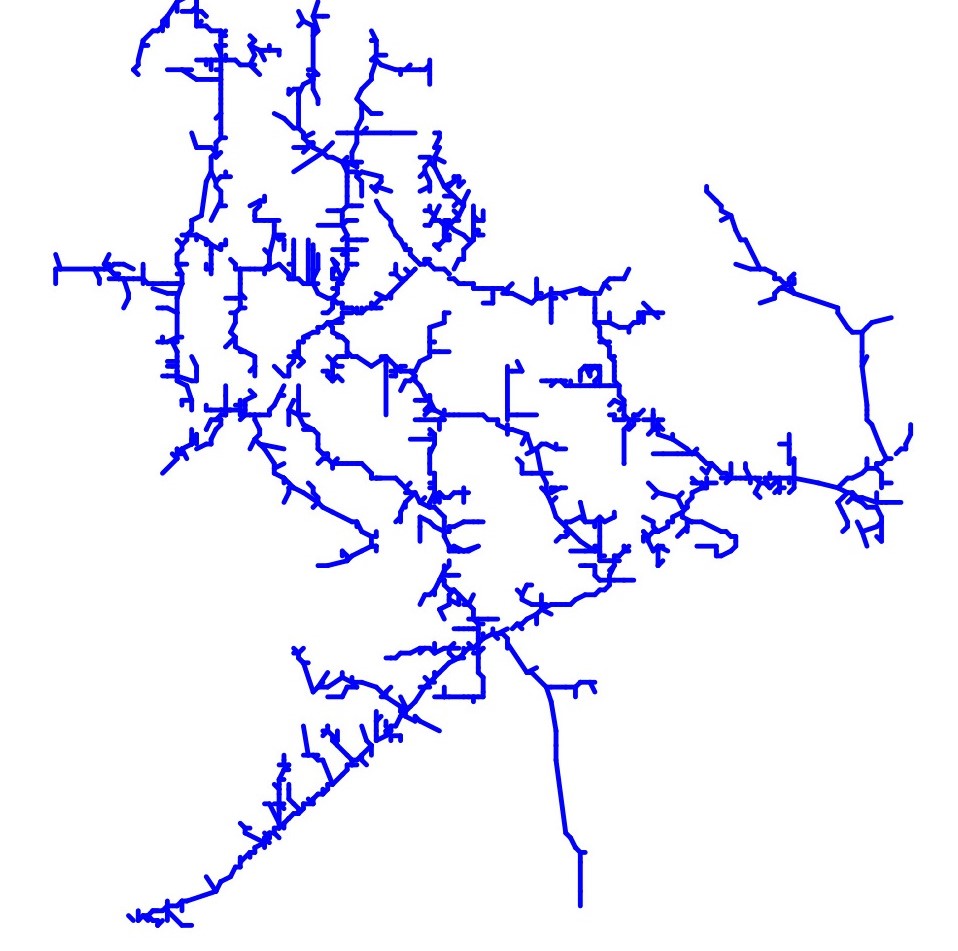}
        \caption{8500-node  feeder in OpenDSS}
        \label{fig:figure1-1}
    \end{subfigure}
    \hfill
    \begin{subfigure}[b]{0.345\textwidth}
        \centering
        \includegraphics[width=\textwidth]{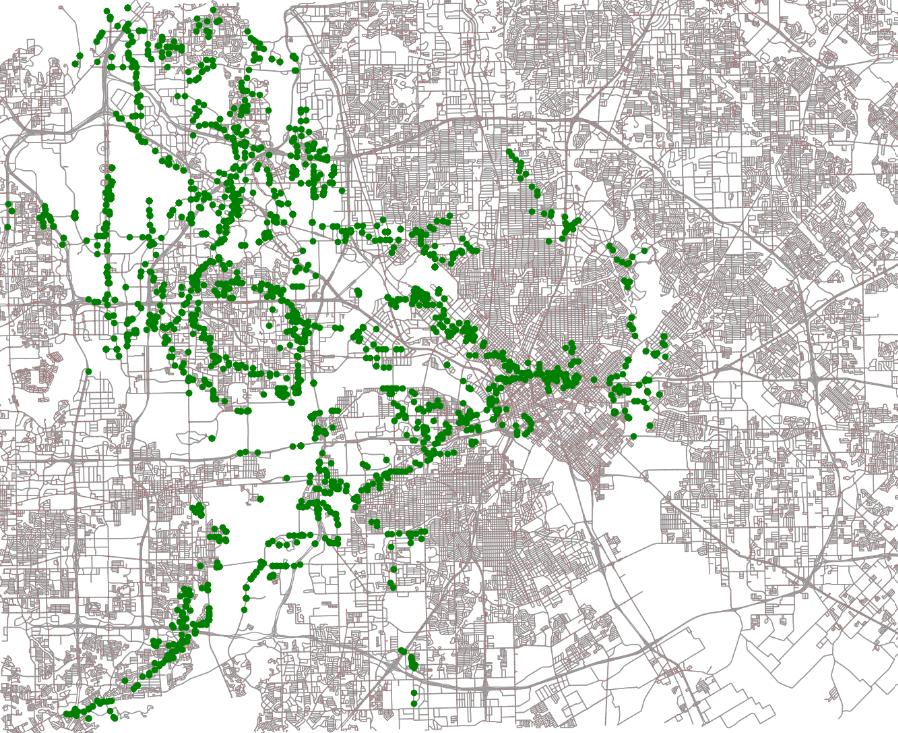}
        \caption{8500 nodes projected on transportation net}
        \label{fig:figure1-3}
    \end{subfigure}
    \hfill
    \begin{subfigure}[b]{0.345\textwidth}
        \centering
        \includegraphics[width=\textwidth]{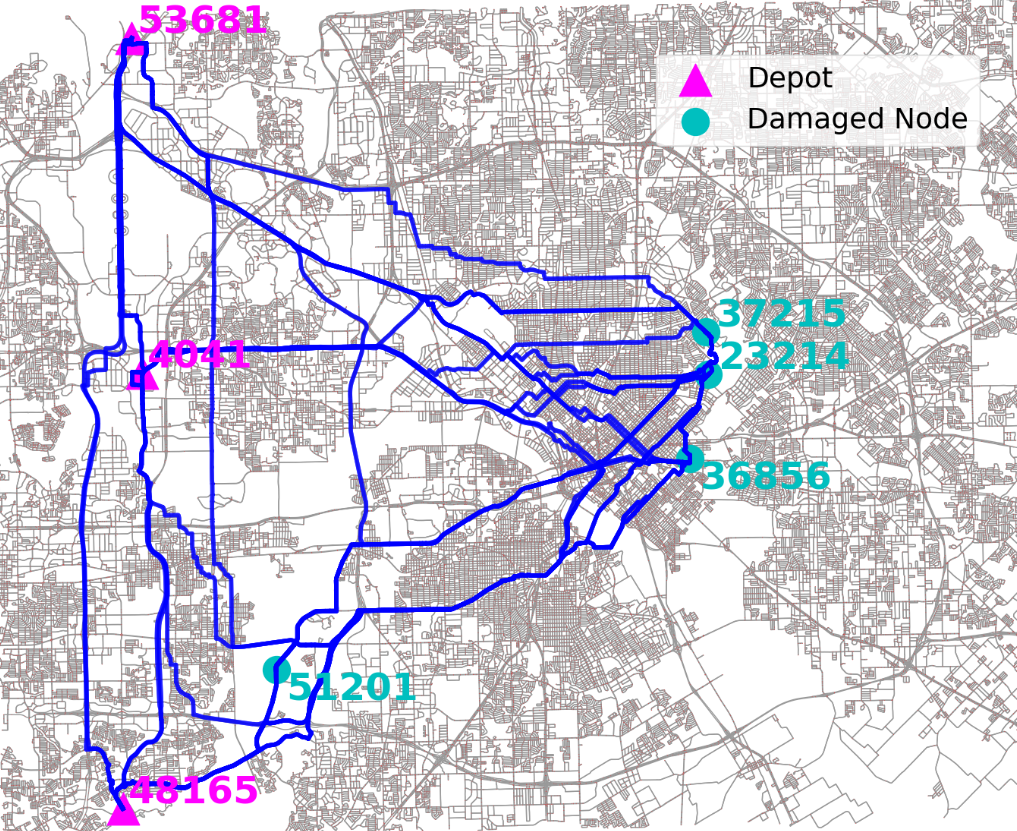}
        \caption{Shortest path among  the damaged nodes and depots.}
        \label{fig:figure_Expt-1_SPT}
    \end{subfigure}
    \caption{Projection of the IEEE 8500-node distribution test network on the DFW transportation network. A complete graph constructed from the damaged nodes and depots.}
\end{figure}

In this study, we employ the IEEE 8500-node distribution test feeder \cite{IEEE8500} as a test case to construct an integrated network for our problem. 
We utilize the network described in the open-source distribution system simulator (OpenDSS), to acquire the representation of the power network and the corresponding coordinates of the buses. This test system is composed of both the primary and secondary levels, of which we consider the more critical primary level of the distribution network in our analysis. As a result, we obtain a network comprising of 2,455 nodes and 2,454 edges, as shown in Figure \ref{fig:figure1-1}. For the transportation network, we opt to utilize the road network of the Dallas-Fort Worth (DFW) area. Leveraging Open Street Maps and NetworkX \cite{BoeingOSMNx}, we procure real road information. We then transform the buses of the primary level of the 8500-node feeder into the spatial framework of the transportation network. 
Damages in the power distribution network include the failure of components such as lines, switches, or transformers in the primary circuit.  When mapping these into the transportation network, the corresponding start buses (nodes in the graph domain) of specific components are modeled as damaged nodes. \hdk{The downstream load $(\mathcal P_i)$ at each of the damaged node is calculated by stimulating the 8500-bus grid using OpenDSS.} \hdk{The green-colored nodes in Figure~\ref{fig:figure1-3} represent the projected power network nodes. The integrated transportation–power network is initially modeled as a large-scale graph consisting of 54,857 nodes and 154,757 edges. To make the problem computationally tractable for the optimization solver, we compute the shortest paths and construct a reduced complete graph over a selected subset of nodes and depots, as shown in Figure~\ref{fig:figure_Expt-1_SPT}.}

\subsection{Simulation and Results}
In this section, we test the performance of our algorithm on a combined network of power and real roads around the Dallas-Fort-Worth (DFW) area. For this study, we have selected only the primary level (critical) of the 8500-bus network, which includes a total of 2455 nodes and  2454 edges. For the transportation network, we utilized the road network around the DFW, obtained using NetworkX and OpenStreetMap \cite{BoeingOSMNx}. We chose this region and required area based on the spatial requirements of the 8500-bus network, which requires a minimum of 2100 square kilometers. Accordingly, we selected a corresponding region within the DFW area.

The next step involves mapping the local nodes from the power network to the transportation network. We converted the local coordinates (from the power network) into geodesic coordinates by adjusting them with specific x and y offsets, based on the bottom-left corner of the transportation network. From NetworkX and OpenStreetMap, we obtained detailed information about graph nodes and edges, including attributes such as OSMIDs, longitudes, latitudes, x and y coordinates of nodes, graph connectivity, and edge lengths. For our analysis, we only focus on the node coordinates (latitude and longitude) and the connectivity of edges.


For our simulations, we used a Windows 11 machine equipped with a 13th Gen Intel Core i7-1365U processor, featuring 10 physical cores and 12 logical threads, and 16 GB of installed RAM. The optimization was performed using Gurobi Optimizer version 11.0.2, with up to 12 threads utilized during execution.

\subsubsection{Case Study 1}

\begin{figure}[t]
    \centering
    \includegraphics[width=0.75\textwidth]{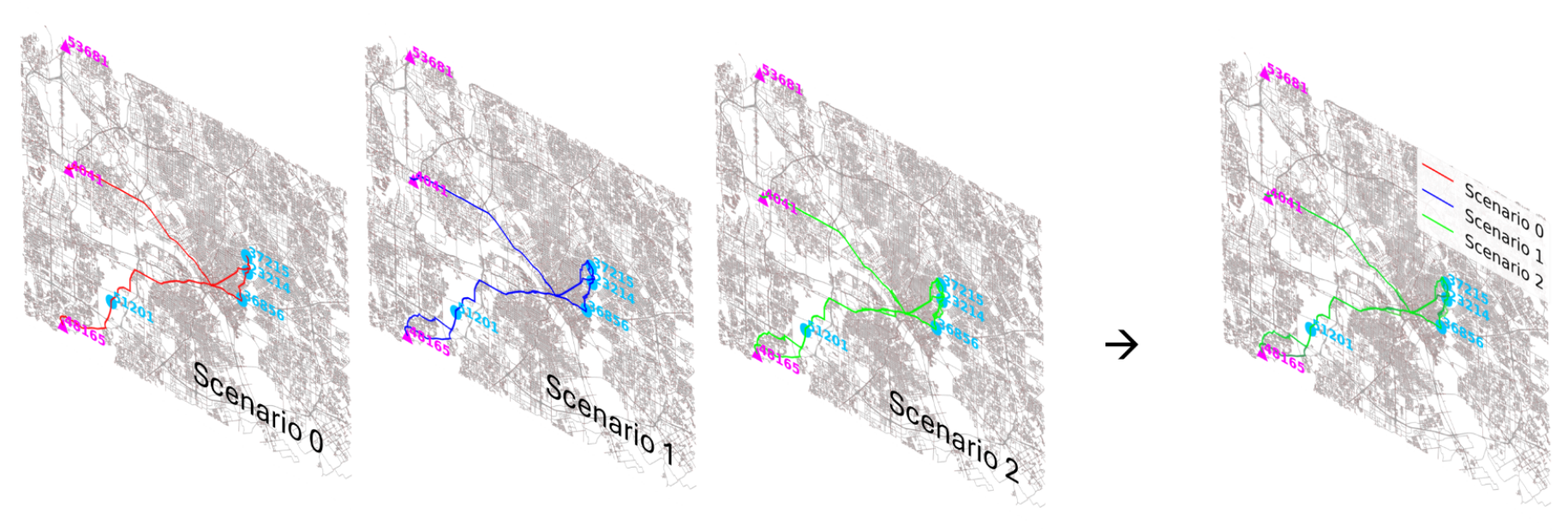}
\caption{Second stage decisions: Crew-2 (Line crew) routing within different scenarios}
    \label{fig:3DCrewRouting}
\end{figure}

\begin{figure}[htpb!]
    \centering
    \begin{subfigure}[b]{0.36\textwidth}
        \centering
        \includegraphics[width=\textwidth]{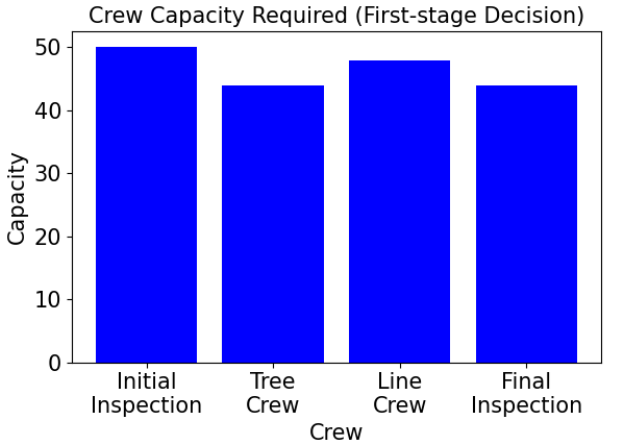}
        \caption{}
    \end{subfigure}
    \hspace{0.9 cm}
    \begin{subfigure}[b]{0.57\textwidth}
        \centering
        \includegraphics[width=\textwidth]{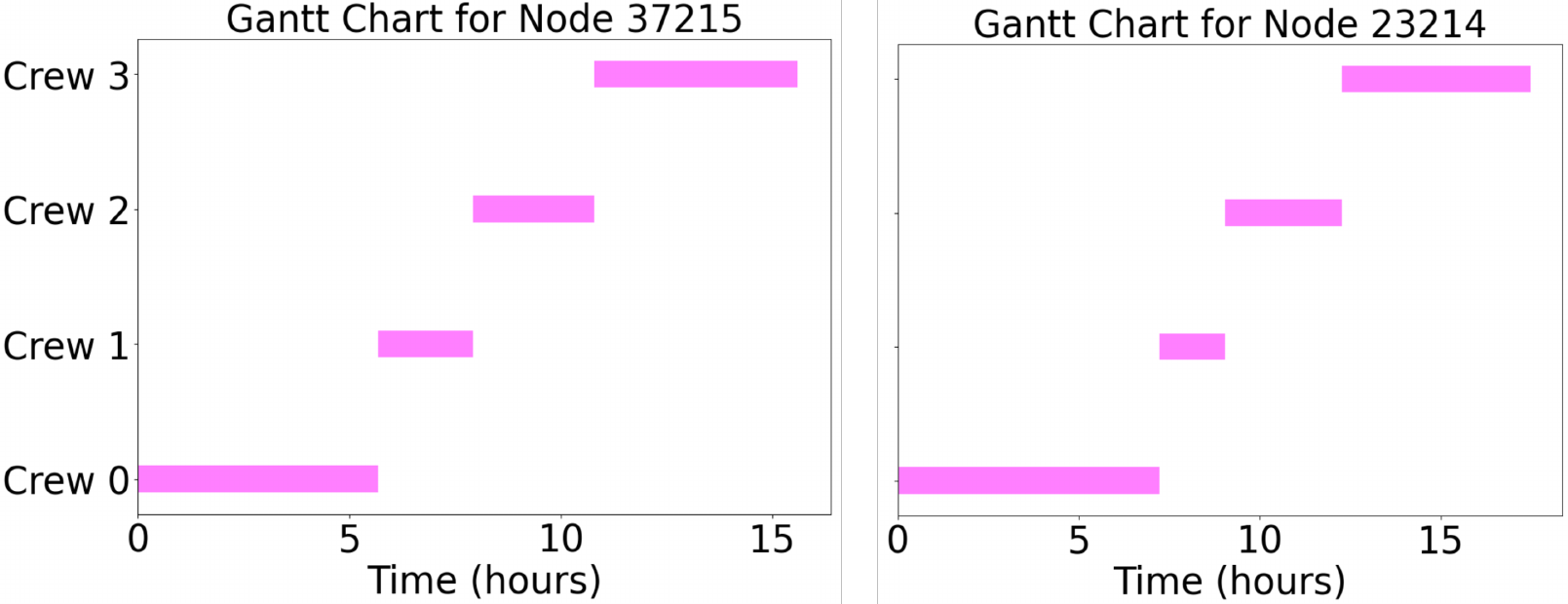}
        \caption{}
    \end{subfigure}\vspace{.66cm}
    \\
    \begin{subfigure}[b]{0.99\textwidth}
        \centering
        \includegraphics[width=\textwidth]{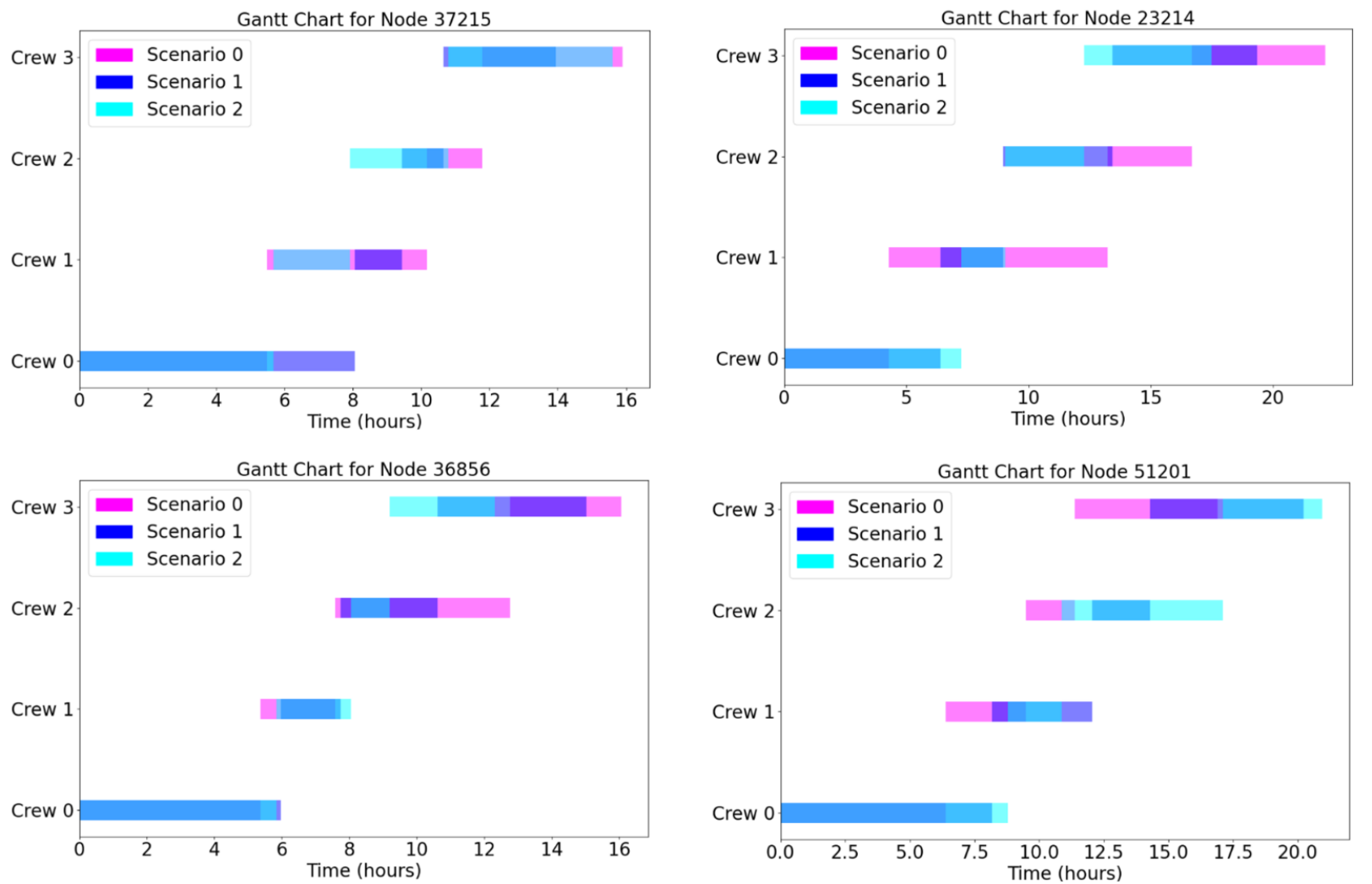}
        \caption{}
    \end{subfigure}\\
    \vspace{0.66cm}
    \begin{subfigure}[b]{0.96\textwidth}
        \centering
        \includegraphics[width=\textwidth]{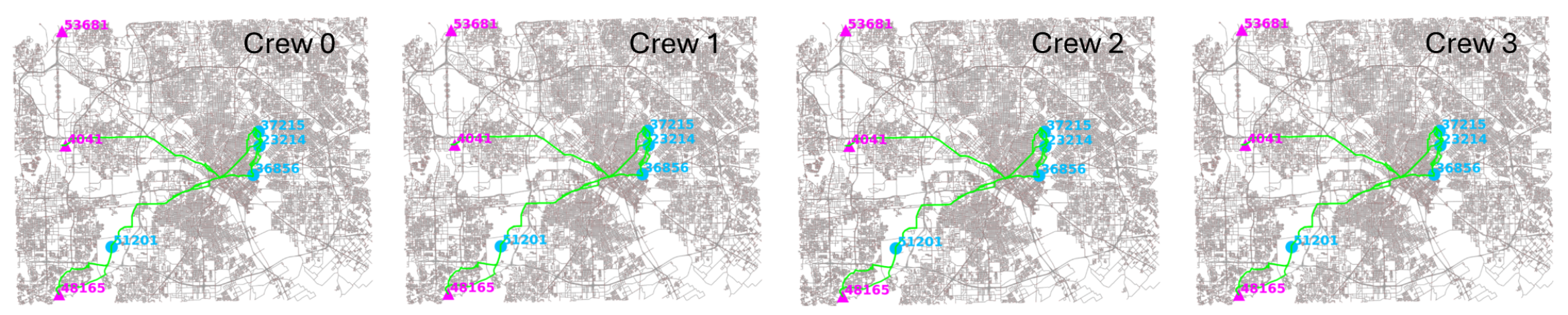}
        \caption{}
    \end{subfigure}
    \caption{\hdk{{\bf Planning Decision After Solving Stochastic Restoration Problem.} {\bf (a)} First-stage decision showing the number of required crews by type. {\bf (b)} Gantt chart illustrating the crew schedule for a single scenario. {\bf (c)}  Gantt chart combining all three scenarios, representing the second-stage scheduling decisions. {\bf (d)} Optimal routing of different crews (Stage 2 decision). The crew assignment and repair sequences are respected. 
    In  plots (b), (c), and (d), crews are labeled as Crew 0 (Initial Inspection), Crew 1 (Tree), Crew 2 (Line), and Crew 3 (Final Inspection).}}
    \label{fig:CaseStudy_1}
\end{figure}


We have considered multiple test cases but for representation purposes, we show the results for three depots and four damaged nodes. In these scenarios, transportation networks are randomly disrupted and the  repair demand and repair times values are randomly selected using the uniforms and lognormal distributions \cite{ZhuZhouYanChen}. Table \ref{tab:scenarios} includes details of how different scenarios are constructed from repair time (in hours) and repair demands. The downstream load $(P_i)$ at each of the damaged node is calculated by stimulating the 8500-bus grid using OpenDSS. Each damaged node is fixed using different crews as shown vertically in Table \ref{tab:scenarios}). For example, Table \ref{tab:scenarios} lists four crews in order to fix node 37215. Each crew has three different repair times and three repair demands corresponding to the three scenarios. Different scenarios generation is discussed in Section \ref{subsec:scenario_gen}. 



\begin{table}[b!]
\centering
\caption{Different scenarios from uncertain repair time and repair demands at the damaged nodes}\label{tab:scenarios}

\begin{tabular}{|c|c|lll|lll|c|}
\hline
\multirow{2}{*}{\textbf{\begin{tabular}[c]{@{}l@{}}Crew\end{tabular}}} & \multicolumn{1}{c|}{}                                                              & \multicolumn{3}{c|}{\textbf{\begin{tabular}[c]{@{}c@{}}Repair time\\ (hours)\end{tabular}}} & \multicolumn{3}{c|}{\textbf{\begin{tabular}[c]{@{}c@{}}Repair \\ demand\end{tabular}}} & \multicolumn{1}{c|}{\multirow{2}{*}{\textbf{\begin{tabular}[c]{@{}c@{}} Power \\ rest-\\ored \\ (kW)\end{tabular}}}} \\ \cline{2-8}
                                                                                         & \textbf{\begin{tabular}{c}
\diagbox[width=5.4em,height=3.6em]{\hspace{-0.36cm}Nodes}{Scenario \hspace*{-0.45cm}\\ (S)} \\
\end{tabular}} & \multicolumn{1}{l|}{\textbf{0}}      & \multicolumn{1}{l|}{\textbf{1}}     & \textbf{2}     & \multicolumn{1}{l|}{\textbf{0}}    & \multicolumn{1}{l|}{\textbf{1}}    & \textbf{2}   & \multicolumn{1}{c|}{}     \\ \hline
\textbf{C0}                                                                              & \multirow{4}{*}{\textbf{37215}}                                                    & \multicolumn{1}{l|}{5.5}             & \multicolumn{1}{l|}{8.1}            & 5.7           & \multicolumn{1}{l|}{17}             & \multicolumn{1}{l|}{16}             & 8            & \multirow{4}{*}{213.6}                                                                                               \\ \cline{1-1} \cline{3-8}
\textbf{C1}                                                                              &                                                                                    & \multicolumn{1}{l|}{4.7}            & \multicolumn{1}{l|}{1.4}           & 2.2          & \multicolumn{1}{l|}{5}            & \multicolumn{1}{l|}{6}             & 8            &                                                                                                                     \\ \cline{1-1} \cline{3-8}
\textbf{C2}                                                                              &                                                                                    & \multicolumn{1}{l|}{1.6}            & \multicolumn{1}{l|}{1.2}           & 2.9           & \multicolumn{1}{l|}{6}             & \multicolumn{1}{l|}{15}             & 12           &                                                                                                                     \\ \cline{1-1} \cline{3-8}
\textbf{C3}                                                                              &                                                                                    & \multicolumn{1}{l|}{4.1}               & \multicolumn{1}{l|}{3.3}            & 4.8            & \multicolumn{1}{l|}{9}             & \multicolumn{1}{l|}{18}             & 6            &                                                                                                                     \\ \hline 
\textbf{C0}                                                                              & \multirow{4}{*}{\textbf{23214}}                                                    & \multicolumn{1}{l|}{4.3}             & \multicolumn{1}{l|}{6.4}            & 7.2           & \multicolumn{1}{l|}{10}            & \multicolumn{1}{l|}{8}             & 10            & \multirow{4}{*}{223.7}                                                                                              \\ \cline{1-1} \cline{3-8}
\textbf{C1}                                                                              &                                                                                    & \multicolumn{1}{l|}{8.9}            & \multicolumn{1}{l|}{2.6}           & 1.8           & \multicolumn{1}{l|}{7}            & \multicolumn{1}{l|}{8}             & 19            &                                                                                                                     \\ \cline{1-1} \cline{3-8}
\textbf{C2}                                                                              &                                                                                    & \multicolumn{1}{l|}{3.5}            & \multicolumn{1}{l|}{4.5}           & 3.2           & \multicolumn{1}{l|}{9}            & \multicolumn{1}{l|}{13}             & 14            &                                                                                                                     \\ \cline{1-1} \cline{3-8}
\textbf{C3}                                                                              &                                                                                    & \multicolumn{1}{l|}{5.5}             & \multicolumn{1}{l|}{5.9}            & 5.2           & \multicolumn{1}{l|}{8}            & \multicolumn{1}{l|}{10}             &     10        &                                                                                                                     \\ \hline 
\textbf{C0}                                                                              & \multirow{4}{*}{\textbf{36856}}                                                    & \multicolumn{1}{l|}{5.4}             & \multicolumn{1}{l|}{6.0}            & 5.8            & \multicolumn{1}{l|}{5}             & \multicolumn{1}{l|}{12}             & 7           & \multirow{4}{*}{287.6}                                                                                              \\ \cline{1-1} \cline{3-8}
\textbf{C1}                                                                              &                                                                                    & \multicolumn{1}{l|}{2.2}            & \multicolumn{1}{l|}{1.8}           & 2.2           & \multicolumn{1}{l|}{8}            & \multicolumn{1}{l|}{18}             & 12            &                                                                                                                     \\ \cline{1-1} \cline{3-8}
\textbf{C2}                                                                              &                                                                                    & \multicolumn{1}{l|}{5.2}            & \multicolumn{1}{l|}{2.9}             & 1.1           & \multicolumn{1}{l|}{14}             & \multicolumn{1}{l|}{6}             & 14            &                                                                                                                     \\ \cline{1-1} \cline{3-8}
\textbf{C3}                                                                              &                                                                                    & \multicolumn{1}{l|}{3.3}             & \multicolumn{1}{l|}{4.4}            &   3.1          & \multicolumn{1}{l|}{15}             & \multicolumn{1}{l|}{10}             &    14         &                                                                                                                     \\ \hline
\textbf{C0}                                                                              & \multirow{4}{*}{\textbf{51201}}                                                    & \multicolumn{1}{l|}{6.4}             & \multicolumn{1}{l|}{8.2}            & 8.8           & \multicolumn{1}{l|}{8}             & \multicolumn{1}{l|}{14}             & 9           & \multirow{4}{*}{8.9}                                                                                              \\ \cline{1-1} \cline{3-8}
\textbf{C1}                                                                              &                                                                                    & \multicolumn{1}{l|}{3.1}            & \multicolumn{1}{l|}{3.9}           & 2.1             & \multicolumn{1}{l|}{13}            & \multicolumn{1}{l|}{8}             & 5            &                                                                                                                     \\ \cline{1-1} \cline{3-8}
\textbf{C2}                                                                              &                                                                                    & \multicolumn{1}{l|}{1.9}            & \multicolumn{1}{l|}{2.2}           & 6.2             & \multicolumn{1}{l|}{15}             & \multicolumn{1}{l|}{6}             & 8            &                                                                                                                     \\ \cline{1-1} \cline{3-8}
\textbf{C3}                                                                              &                                                                                    & \multicolumn{1}{l|}{5.5}             & \multicolumn{1}{l|}{5.9}            &  3.8           & \multicolumn{1}{l|}{12}             & \multicolumn{1}{l|}{5}             & 8            &                                                                                                                     \\ \hline
\end{tabular}
\end{table}

After determining the set of scenarios and finalizing all input parameters, we solve the stochastic optimization problem following Algorithm \ref{algorithm1}. The first-stage decisions are presented in Figure \ref{fig:CaseStudy_1}(a), which shows the required crew capacities across all potential scenarios. The second-stage decisions are illustrated in Figure \ref{fig:CaseStudy_1}(b), (c), and (d) using Gantt charts and optimal routing plans for each crew type. Figure \ref{fig:CaseStudy_1}(b) provides the Gantt chart for a representative scenario, while Figure \ref{fig:CaseStudy_1}(c) combines the schedules across all scenarios. The optimal routing paths for individual crews, such as repair crew 2 (line crew), are also depicted across scenarios. These routing plans serve as actionable guides once a specific scenario is realized.

To ensure a coordinated and efficient repair process, we generate detailed time schedules for each crew at all damaged nodes, as shown in Figures \ref{fig:CaseStudy_1}(b) and (c). These schedules are prepared for all considered scenarios. For example, at node 37215 under scenario 0 (magenta), the inspection crew begins work at hour 0 and completes by approximately hour 5.5. This is followed by the tree trimming crew for 4.6 hours and the line crew for 1.6 hours. Finally, the final inspection crew is scheduled to start around hour 11.8, assuming full crew availability and no delays. These comprehensive and synchronized repair schedules across all damaged nodes represent the core output of the second-stage decisions.

\subsubsection{Case Study 2}

\begin{figure}[htpb!]
    \centering
    \begin{subfigure}[b]{0.99\textwidth}
        \centering
        \includegraphics[width=\textwidth]{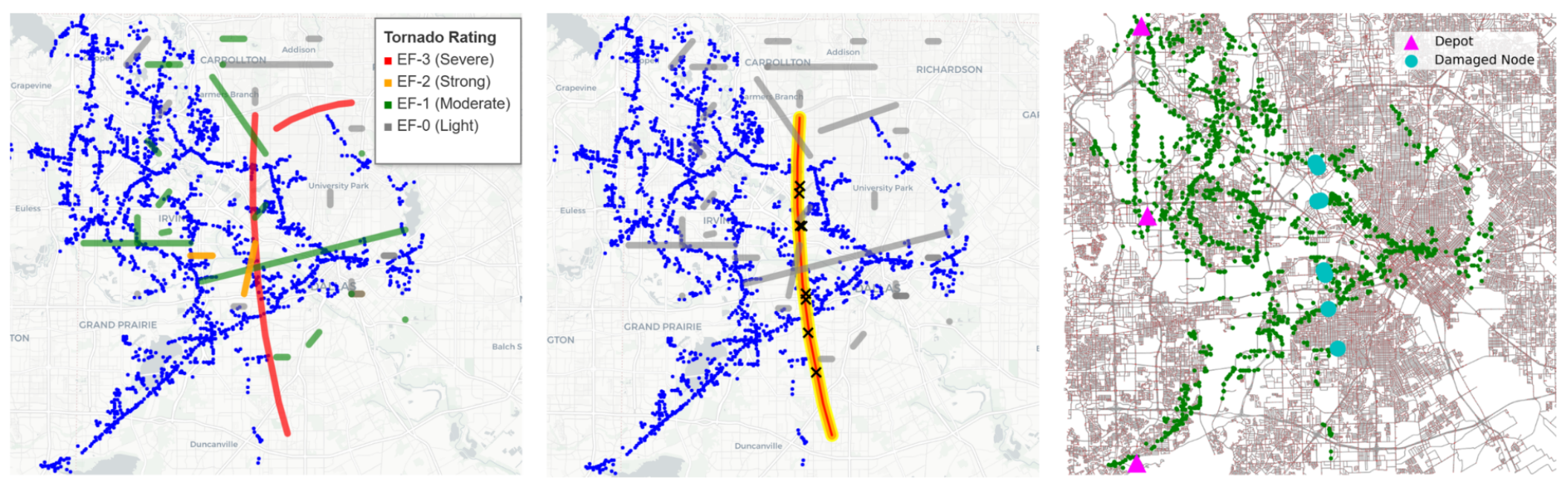}
        \caption{}
    \end{subfigure}\\
    \vspace{.66cm}
    \begin{subfigure}[b]{0.3\textwidth}
        \centering
        \includegraphics[width=\textwidth]{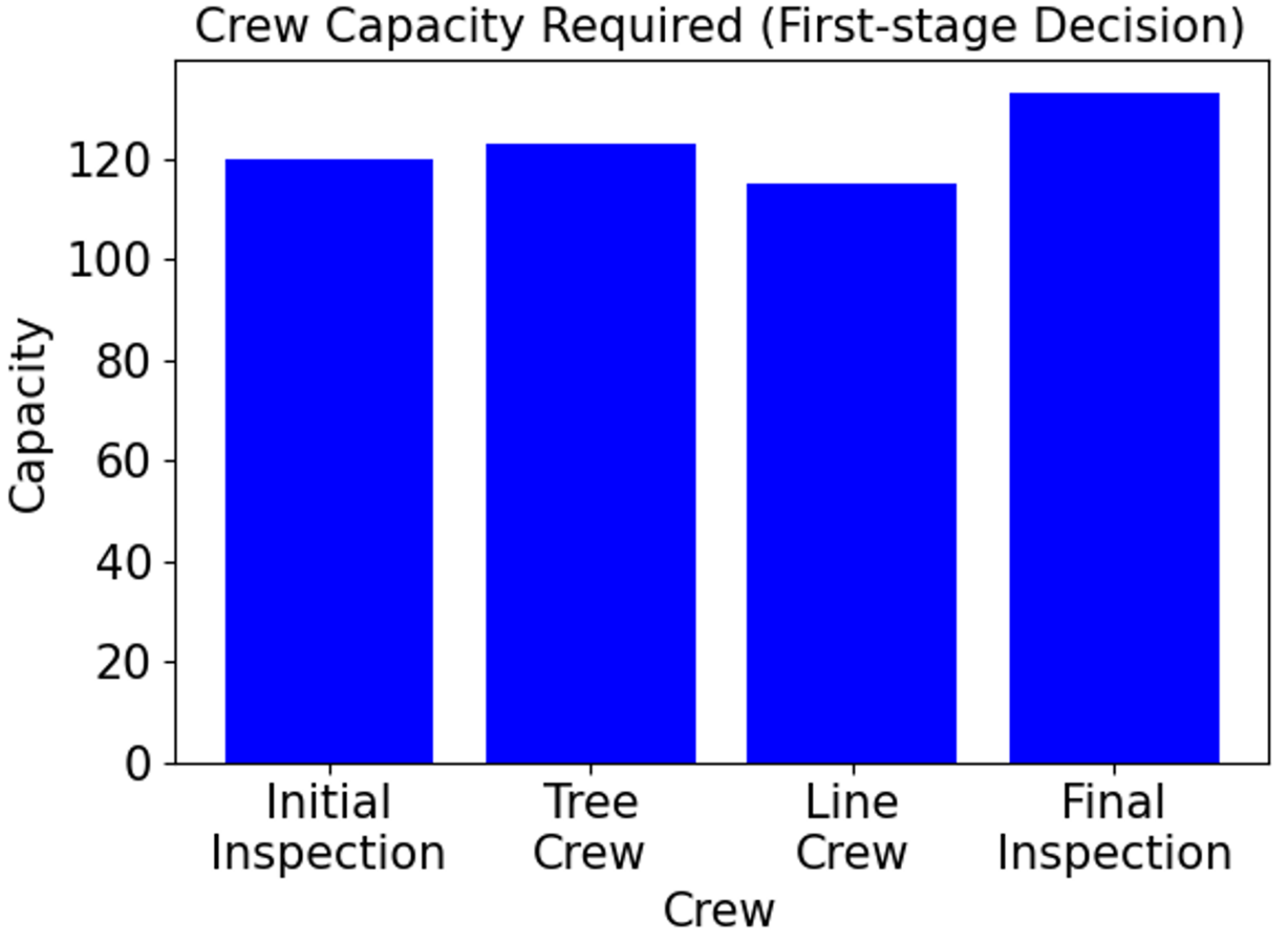}
        \caption{}
    \end{subfigure}
    \hspace{0.9 cm}
    \begin{subfigure}[b]{0.63\textwidth}
        \centering
        \includegraphics[width=\textwidth]{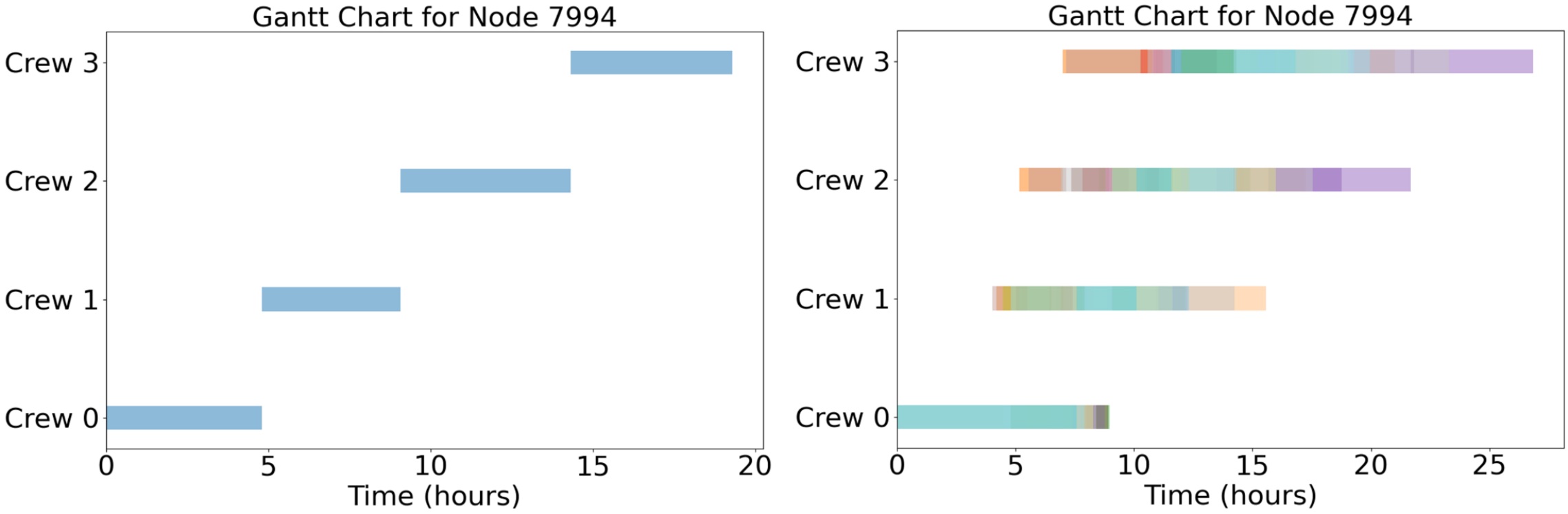}
        \caption{}
    \end{subfigure}
    \caption{\hdk{{\bf  Planning Decisions Based on Actual Tornado Events in the Dallas Area Using Algorithm~\ref{algorithm1}.} {\bf (a)} Tornado events across the DFW area are analyzed based on severity. The most severe event is selected, and its impacted power grid nodes are identified. The third panel illustrates these damaged nodes and available depots on the integrated power grid and transportation network. {\bf (b)} First-stage decision showing the required number of crews by type, determined based on disaster intensity and considering all the scenarios. {\bf (c)}   Second-stage scheduling decisions visualized using Gantt charts. The first chart represents a single scenario, while the second chart combines all scenarios to show aggregate scheduling decisions. In  plots (c), crews are labeled as Crew 0 (Initial Inspection), Crew 1 (Tree), Crew 2 (Line), and Crew 3 (Final Inspection).}}
    \label{fig:CaseStudy_2}
\end{figure}

\hdk{The second case study is developed using real tornado data within the DFW area. This dataset, obtained from the National Weather Service \cite{nws_tornado_climatology_2025}, contains detailed information about each tornado event, including its intensity rating (Enhanced Fujita scale), geographic path (start and end latitude-longitude), duration, and associated damages. As shown in Figure \ref{fig:CaseStudy_2}(a), we first visualize all historical tornado events across the region, categorized by severity. Among them, we select a particularly destructive EF-2 tornado for further analysis due to its significant impact on the local infrastructure. We then identify the power grid nodes that fall within the affected zone of this tornado and mark them as damaged. These damaged nodes, along with available depots, are mapped onto the integrated power distribution and transportation network (Figure \ref{fig:CaseStudy_2}(a), rightmost panel). This spatial setup serves as the input to the optimization model.}

\hdk{To solve the stochastic restoration problem using the previously constructed damaged nodes, depots, and coupled power-transportation network, we apply the proposed Algorithm \ref{algorithm1}. The algorithm determines both first-stage and second-stage decisions under uncertainty. The first-stage decisions, shown in Figure \ref{fig:CaseStudy_2}(b), capture the required capacity (number of crews) for each crew type of crew: inspection, tree removal, line repair, and final inspection, while considering all the potential future scenarios. The second-stage decisions are illustrated in Figure \ref{fig:CaseStudy_2}(c) through Gantt charts. These charts provide detailed crew schedules across time and nodes. The left panel in Figure \ref{fig:CaseStudy_2}(c) shows the crew dispatch plan for a single representative scenario, while the right panel aggregates scheduling outcomes across all 30 scenarios. These scenarios are generated by introducing uncertainty into the transportation network availability, repair time variability, and fluctuating restoration demands. Together, this case study demonstrates how real-world disaster data can be integrated with our optimization framework to develop robust and equitable power grid restoration strategies.}

\section{DISCUSSION}\label{sec:discussion}
\hdk{The proposed multi-stage stochastic optimization framework addresses the complex and time-sensitive challenge of post-disaster power grid restoration by integrating uncertainty in both infrastructure damage and recovery logistics. Our approach explicitly considers realistic sources of uncertainty, including variable repair times, fluctuating repair demands, and disruptions in the transportation network, by generating a diverse set of scenarios informed by actual tornado data, power grid infrastructure, and road network topologies.}

\hdk{Unlike many existing approaches that treat infrastructure systems in isolation, we present a coupled model that jointly considers the power distribution network and the transportation network used by utility crews. To the best of our knowledge, this is one of the first implementations of a multi-agent capacitated vehicle routing formulation over an interdependent infrastructure network for restoration planning, building upon our earlier work \cite{KPEC_HDK}, \cite{Nathan_HDK2025}. This coupling is essential, as disasters such as tornadoes often damage both electricity infrastructure and roadways simultaneously, impacting crew mobility and repair feasibility.}

\hdk{In our framework, scenario generation is grounded in physical data, such as historical tornado paths, network topology, and component-level failure probabilities. Although this process can be further enriched by incorporating additional dimensions such as building criticality, land-use characteristics, or real-time meteorological forecasts, our current methodology establishes a practical and data-driven foundation for generating representative and actionable damage scenarios. The framework employs a two-stage decision process: in Stage 1, we determine the required number and types of crews based on aggregated scenario data, while in Stage 2, we generate optimal crew schedules and routing decisions for each scenario. The use of scenario decomposition and multi-agent Capacitated VRP enables us to manage computational complexity, ensure convergence, and maintain repair sequence feasibility.}

\hdk{In Algorithm \ref{algorithm1}, our decomposition strategy was critical to handling the complexity of solving the assignment and routing subproblems jointly. Without this, convergence to optimal solutions was not attainable given the combinatorial scale of the problem. While we currently assume only four non-overlapping crew types, our modular setup allows for future extension to more diverse crew hierarchies, partial overlaps, or even dynamic task reassignment during restoration execution.}

\hdk{To validate the applicability of our model, we conducted two case studies. In the first, we considered three representative scenarios where damaged nodes and corresponding failure attributes were randomly generated using lognormal and uniform distributions. This limited set was chosen to facilitate clear discussion and visualization of simulation results. In the second case study, we leveraged real tornado event data from the DFW region to construct realistic failure scenarios and evaluate the model’s performance in a more practical restoration setting. The tornado dataset included spatial and severity attributes, allowing us to identify damaged grid nodes and disrupted road segments. These were used to construct relevant damage scenarios, which were then solved using our proposed framework. Although more sophisticated scenario generation, potentially involving domain experts from meteorology, civil infrastructure, or emergency response, can enhance realism, our focus here was on demonstrating the ability of our optimization model to produce fair, feasible, and robust restoration schedules.}



\section{METHODS}\label{sec:methods}
In this section, we present the proposed methodology in detail. We begin by describing the mathematical formulation, followed by a discussion of the two-stage structure and the overall approach used to address the problem. Finally, we provide an outline of the algorithm used to implement the solution framework.



\subsection{Mathematical formulation} \label{subsec:math_form}

We address the repair and restoration problem  over a  combined transportation and power network graph $\mathrm{G(\mathbb{V},\mathbb{E})}$, formed by a node set $\mathbb{V}$, which includes both damaged nodes $\mathbb{N}$ and depots $\mathbb{D}$, and an edge set  $\mathbb{E}$ representing the connections between these nodes utilizing the available road network. This includes undirected edges that represent the actual road network. 
Specifically, we focus on the subset of nodes $\mathbb{N}$ that are damaged and require restoration. This undirected graph allows us to model the network's structure and the interconnected components effectively. The recovery process for any damaged node involves a sequence of repair crews (total set of repair crews: $\mathbb{K}$), where each crew is tasked with different types of repairs and the specific requirements of the damaged nodes. \hdk{The type of crew $k\in \mathbb{K}$ is assigned to each damaged node by ensuring the restoration process is efficient and that the appropriate resources are utilized for each task with the least idle time for any crew.}

\hdk{The repair demand $\mathcal{D}$ at different damaged nodes is the number of humans required for the repair work. This capacity is dictated by the severity of the damage, ensuring that sufficient manpower is available to address all restoration needs adequately. At each damaged node, 
$\mathcal T$ denotes the time required (in hours) to complete the repair. To capture the inherent uncertainties in the network restoration process, we construct a set of scenarios denoted by  $\mathcal{S}$. Each scenario in $\mathcal{S}$ reflects variations in node-specific repair demand, repair time requirements, and the availability of road networks. Given that road infrastructure can also be compromised during a storm, we assume that the set of traversable roads differs across scenarios, thereby altering feasible travel paths between nodes.}

\subsubsection{First Stage: Strategic Allocation and Assignment of Repair Crews}
This problem involves assignment and sequencing optimization, with the primary objective of the first-stage optimization being to determine the optimal allocation and sequencing of repair crews ($k \in \mathbb{K}$) for restoring damaged nodes $\mathbb{N}$. The model uses as inputs the estimated repair times ($\mathcal{T}$ in hours), the power restoration potential ($\mathcal{P}$ in kW) at each damaged node, and the repair demand ($\mathcal{D}$) to finalize strategic crew assignment decisions.

\noindent{\hdk{Objective function: }} The objective function in the first stage aims to maximize power restoration while minimizing total repair time across the power grid. Specifically, it incorporates incentives for restoring nodes with higher power potential ($\mathcal{P}$) and prioritizing those with lower repair times ($\mathcal{T}$). The primary decision variable is the integer variable representing the total number of crew members to be mobilized, denoted by $x_k$ for each $k \in \mathbb{K}$. Additionally, the assignment of a specific type and amount of crew to each damaged node under scenario $s \in \mathcal{S}$, is capture by using another continuous  variable $y_{s,i,k}$ for node $ i \in \mathbb N$, and crew $k \in \mathbb{K}$. The objective for the first stage is formalized as,
\begin{align*}
    \text{minimize} \quad  \mathrm c \ \sum_{k\in\mathbb{K}} \textrm{C}^1_kx_{k} - \frac{1}{|\mathcal{S}|}\left\{\sum_{s\in \mathcal{S}}\sum_{i\in \mathbb N}\sum_{k\in \mathbb{K}}\mathcal{P}_{i}  \ y_{s,i,k}\nonumber - \sum_{s\in \mathcal{S}}\sum_{i\in \mathbb N}\sum_{k\in \mathbb{K}}\mathcal{T}_{s,i,k}  \ y_{s,i,k}\right\}, 
\end{align*}    
where $\mathrm c $ is an appropriate  scalar multiplier, chosen to be sufficiently large to balance the trade-off between maximizing power restoration, minimizing repair time, and  the total number of crew members within the combined objective function. \hdk{$\textrm{C}^1_k$ represents the per-hour labor cost associated with crew type $k$. For the Dallas, Texas region, typical rates chosen as: tree-trimming crews charge approximately \$60-70 per person per hour, line repair crews range from \$50 to \$100 per hour, and both initial and final inspection crews, often comprising specialized technicians or engineers, command higher rates between \$175 and \$225 per person per hour.}


\hdk{Constraint set:} The first set of constraints ensures that the total available crew capacity is sufficient to meet the repair demand requirements at all damaged nodes.
\begin{align*}
    \sum_{i\in \mathbb N} y_{s,i,k} \leq x_k \quad \text{for all }  s\in \mathcal{S}, \ k \in \mathbb{K}.
\end{align*}
The second set of constraints enforces that, for each scenario and crew type, the assigned crew resources are adequate to fulfill the repair demands at all damaged nodes, and is mathematically expressed as,
\begin{align*}
    \mathcal{D}_{s,i,k} \leq y_{s,i,k} \quad \text{for all }  s\in \mathcal{S}, \  i\in \mathbb{N}, \ k \in \mathbb{K}.
\end{align*}
Lastly, the non-negativity constraints on both the set of variables $x_k$ and $y_{s,i,k}$ for $k \in \mathbb{K}, i\in\mathbb{N}, $ and $s\in \mathcal{S}$ are enforced. 

\subsubsection{\hdk{Second Stage: Capacitated Vehicle Routing for Post-Disaster Crew Dispatch}} In the second stage, the model focuses on the operational decisions required under specific scenarios. The primary focus is on route planning for crew deployment, captured through a binary decision variable  $z_{s,i,j,k}$  indicates whether crew $k\in \mathbb{K}$ moves from node $i$ to node $j$ within scenario $s\in \mathcal{S}$.  

\noindent\hdk{Objective function:} The objective is to construct efficient crew routes that minimize total travel cost while ensuring that the repair requirements at the damaged nodes are satisfied. These routing decisions follow the repair sequence determined in stage 1, ensuring consistency between strategic crew assignment and tactical deployment across scenarios. The objective function for this stage is given as,

\begin{align*}
   \text{minimize}\ \frac{1}{|\mathcal{S}|} \ \sum_{s\in \mathcal{S}}\sum_{(i, j)\in \mathbb E}\sum_{k\in \mathbb{K}}\ \mathrm{C}_{s,i,j,k}^2 \ z_{s,i,j,k},
\end{align*}
where $\mathrm{C}_{s,i,j,k}^2$ is the travel cost for crew $k$ from location $i$ to $j$, within scenario $s\in \mathcal{S}.$

\hdk{Constraint set:} The first set of constraints makes sure that each damaged node  $i \in \mathbb{N}$ is visited once as long as there is a repair demand $\mathcal{D}_{s,i,k} $ for crew $k$ within scenario $s$. 
\begin{align*}
    &\sum_{j\in \mathbb{N}} z_{s,i,j,k} = 1 \quad  \text{for } s\in\mathcal{S},\  k\in \mathbb{K}, \ i \in \mathbb{N}\ \textrm{ if } \ \mathcal{D}_{s,i,k} > 0,\\
    &\sum_{i\in \mathbb{N}} z_{s,i,j,k} = 1  \quad \text{ for } s\in\mathcal{S},\  k\in \mathbb{K},\  j \in \mathbb{N}\ \textrm{ if } \ \mathcal{D}_{s,i,k} > 0.
\end{align*}
The second set of constraints makes sure that each crew $k$ starts and ends at one of the depots in $\mathbb{D}$.
\begin{align*}
    \sum_{j\in \mathbb{N}} z_{s,i,j,k} = 1 \quad &\text{for } s\in\mathcal{S}, \ k\in \mathbb{K}, \ i \in \mathbb{D},\\
    \sum_{i\in \mathbb{N}} z_{s,i,j,k} = 1  \quad &\text{for } s\in\mathcal{S},\  k\in \mathbb{K},\  j \in \mathbb{D}.
\end{align*}
The third set of constraints establishes flow conservation within the network. It specifies that once a repair crew $k$ arrives at a damaged component, they must proceed to the next location after completing the repairs, ensuring a continuous flow of operations.
\begin{align*}
        &\sum_{(i,j) \in \mathbb E}z_{s,i,j,k} = \sum_{(j,i) \in \mathbb  E}z_{s,j,i,k} \quad \text{ for } j\in \mathbb{N},\  s \in \mathcal{S}, \ k \in \mathbb{K}.
\end{align*}
The fourth set of constraints is the Miller-Tucker-Zemlin (MTZ) subtour elimination constraint \cite{MTZ}. 
\begin{align*}
    u_{s,i,k} - u_{s,j,k} + |\mathbb{N}| \ z_{s,i,j,k} \leq |\mathbb{N}| -1 \quad \textrm{ for } s\in \mathcal{S}, \ i\in\mathbb{N},\  j\in\mathbb{N}, \ (i,j) \in \mathbb  E, \textrm{ and } k \in \mathbb{K},
\end{align*}
where variable $u_{s,i,k}$ represents the position of  node $i$ in the route for crew $k$ in scenario $s$. This constraint ensures that, for each scenario, the model avoids disconnected cycles (subtours) by enforcing a consistent order of visits to nodes in the routing path.

 Finally, the last set of constraints enforce non-negativity requirements and ensure that decision variables are binary. 
\begin{align*}
    &u_{s,i,k} \geq 0 && \hspace{-3.96cm} \text{ for } i \in \mathbb{N}, s \in \mathcal{S}, k \in \mathbb{K},\\
    &z_{s,i,j,k} \in \{0,1\} && \hspace{-3.96cm} \text{ for } s \in \mathcal{S}, (i,j)\in \mathbb E, k \in \mathbb{K}.
\end{align*}

\hdk{We emphasize that the proposed two-stage formulation introduces significantly greater complexity compared to traditional approaches. Unlike much of the existing literature, which typically addresses a single source of uncertainty due to computational challenges, our model incorporates multiple, interdependent uncertainties, including variable repair demands, repair times, and dynamic road availability. Furthermore, to the best of our knowledge, this is the first extension of grid repair and restoration methods \cite{ArifWangChenWang1, ArifWangWangChen2018, ArifMaWangWangRyanChen2018} to a stochastic graph-based framework that explicitly integrates both transportation and power network.}

\subsection{Two-stage Structure of the Algorithm Proposed} 
\label{subsec:two_stage_structure}

\begin{figure}[htpb!]
     \centering
    \includegraphics[width=0.99\textwidth]{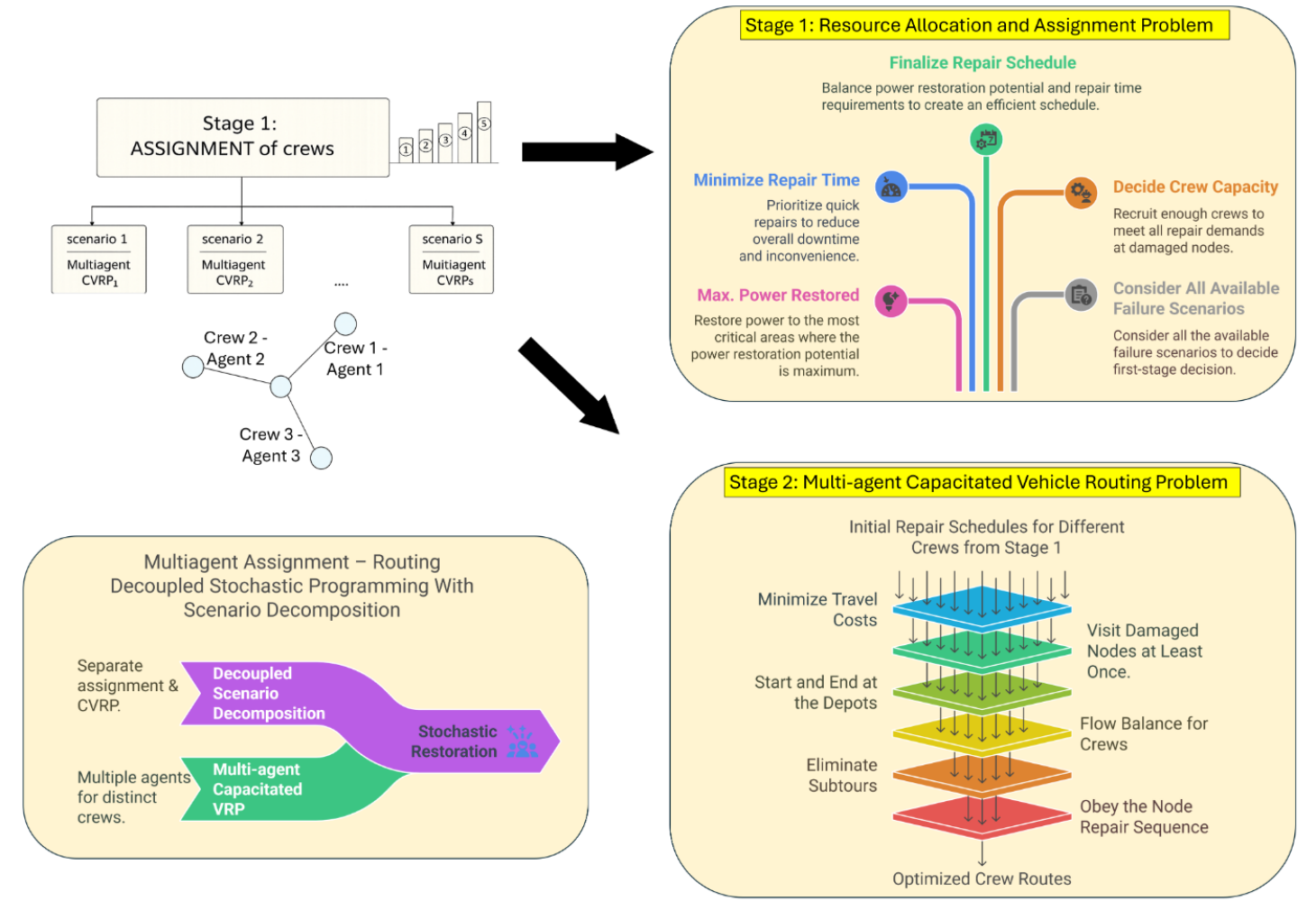}
     \caption{\hdk{Multiagent Assignment-Routing Decoupled Stochastic Programming Framework with Scenario Decomposition. The two-stage optimization model addresses power grid restoration by first assigning repair crews across pre-generated damage scenarios (Stage 1), maximizing power restoration and minimizing repair times. In Stage 2, each scenario is solved using a Multi-Agent Capacitated Vehicle Routing Problem, where different crew types are routed efficiently through the transportation network. The scenario decomposition strategy decouples assignment and routing decisions to enhance scalability and convergence. Distinct agents are assigned to individual repair crews, enabling decentralized coordination and efficient task allocation across the network.}}
    \label{fig:two_stage_structure}
\end{figure}

We address the stochastic grid repair and restoration problem by proactively accounting for a diverse set of potential damage scenarios (as described in Section \ref{subsec:scenario_gen}). These scenarios are generated in advance and represent realistic outcomes of extreme events. Our goal is to pre-position repair crews in a way that enables a rapid and effective response regardless of which scenario eventually materializes.

The decision-making framework is structured as a two-stage stochastic optimization model. In the first stage, we perform a strategic allocation and assignment of repair crews before any scenario occurs. This includes determining the number and type of crews ($k \in \mathbb{K}$) to be mobilized and their assignment to damaged nodes $\mathbb{N}$ based on their restoration potential ($\mathcal{P}$), estimated repair times ($\mathcal{T}$), and repair demands ($\mathcal{D}$). These decisions are scenario-independent and must ensure adequate capacity to address the full range of failure conditions. The objective of this stage is to maximize power restoration and minimize total repair time, while respecting crew capacity and coverage constraints.

In the second stage, once a specific damage scenario is realized, the problem becomes a multi-agent capacitated vehicle routing problem. Each agent corresponds to a distinct type of repair crew (e.g., line crews, tree crews), and the routing decisions ensure that the crews travel from depots to the assigned damaged locations in the most cost-effective and timely manner. This stage respects the repair assignments made in stage 1 and incorporates classical vehicle routing problem constraints such as flow balance, depot return, and subtour elimination (e.g., MTZ constraints \cite{MTZ}). The routing paths are optimized within the constraints of the transportation network, ensuring that each damaged node is visited and serviced, with the goal of minimizing travel costs and expediting restoration. To effectively solve this complex problem, we employ two key strategies:
\begin{itemize}[leftmargin=0.36cm]
    \item \hdk{Scenario Decomposition:  The resource allocation and assignment (Stage 1) and routing (Stage 2) problems are decoupled and solved separately. This decomposition approach was adopted due to convergence challenges encountered when attempting a unified model that combines assignment and the capacitated vehicle routing problem formulations.}
    \item \hdk{Multi-Agent Capacitated Vehicle Routing Problem: Different types of crews are modeled as separate agents with distinct capabilities and resource constraints. This approach allows for realistic modeling of heterogeneous crews, enabling decentralized coordination and efficient task allocation across the network.}
\end{itemize}

\subsection{Algorithm Outline}\label{subsec:algo_outline}
In this section, we explain the structure of our algorithm to address the two-stage stochastic grid restoration problem. 

\begin{algorithm}[H]
\caption{Two-Stage Stochastic Programming With Scenario Decomposition}
\label{algorithm1}
\textbf{Discuss and finalize:} A set of failure scenarios $\mathcal{S}$ (Section~\ref{subsec:scenario_gen}): 
\begin{align*}
    &&\mathcal{D}: &\text{ Set of demands at all the damaged nodes. }\\
    &&\mathcal{T}: &\text{ Repair times at the damaged nodes (in hours). }\hspace{7.2cm}\\
    &&\mathcal{P}: &\text{ Power restoration potential at the damaged nodes (in kW).}
\end{align*}
\textbf{Generate:} Coupled power and transportation network $\mathrm G (\mathbb{V}, \mathbb{E}) $(Section~\ref{subsec:integrated_network}).\\\\
\textbf{Parameters for two-stage grid restoration problem:}
\begin{align*}
&\mathrm{C}_k^1&&: \text{Cost parameter 1 in first stage objective function.}\\
&\mathrm{C}_{s,i,j}^2&&: \text{Cost parameter 2 in second stage objective function.}\hspace{3.6cm}\vspace*{3cm}\\
& \mathcal{S} &&: \text{Set of all scenarios.}\\
&\mathcal{P}_i&&: \text{Downstream load (kW) at node $i$ (from OpenDSS).}\\
&\mathcal{T}_{s,i,k}&&: \text{Repair time at node $i$ by crew $k$ in scenario $s$.}\\
&\mathcal{D}_{s,i,k}&&: \text{Repair demand at node $i$ by crew $k$ in scenario $s$.}\\
&\mathbb{K} &&: \text{Set of crew types.}\\
&\mathbb{N} &&: \text{Set of damaged nodes.}\\
&\mathbb V \triangleq \mathbb{N} \cup \mathbb{D} &&: \text{Set of all nodes.}\\
&\mathbb{E} &&: \text{Set of all edges.}\\
&\mathrm{G(\mathbb V,\mathbb E)} &&: \text{Coupled power and transportation network.}
\end{align*}
\textbf{Decision Variables:}
\begin{align*}
    &x_k&&: \text{($1^{\text{st}}$ stage) Crew $k \in \mathbb{K} $ capacity (integer variable).}\\
&y_{s,i,k}&&: \text{($1^{\text{st}}$ stage) Crew $k $ members assigned to node $i$ within $s$ (continuous variable).}\hspace{-1.2cm}\\
&z_{s,i,j,k}&&: \text{($2^{\text{nd}}$ stage) Binary variable indicating if crew $k$ travels from $i$ to $j$ in scenario $s$.}
\end{align*}

\end{algorithm}

\hdk{The proposed two-stage optimization framework for grid restoration is designed to address uncertainties in repair times ($\mathcal{T}$), repair demands ($\mathcal{D}$), and transportation network disruptions. Algorithm \ref{algorithm1} implements a two-stage multiagent stochastic programming approach with scenario decomposition, where the first stage determines pre-disaster crew allocation across all possible damage scenarios, and the second stage solves a scenario-specific multiagent vehicle routing problem to dispatch crews efficiently once a scenario materializes. Initially, all disruption scenarios $(\mathcal{S})$ are generated. For each scenario, the transportation graph $(\mathrm{G})$ is constructed, downstream loads $(\mathcal{P})$ are simulated using the 8500-bus network in OpenDSS, and repair parameters $(\mathcal{T}, \mathcal{D})$ are sampled from appropriate distributions \cite{ZhuZhouYanChen}. Next, within the first stage, we decide for resource allocation and determine the necessary capacity for each type of crew involved in the restoration process by considering all the damaged scenarios.} Once the resources have been allocated and the required capacities have been determined, the second stage focuses on the scheduling and optimal routing of the crews. The respective decision variables and input parameters are summarized in Algorithm \ref{algorithm1}. This two-stage structure effectively allows for a balance between preparedness and adaptability. 

\section{CONCLUSION}\label{sec:conclude}
\hdk{In this work, we presented a two-stage stochastic optimization framework for power grid restoration that jointly models interdependencies between the power distribution and transportation networks while explicitly accounting for uncertainty in disaster impacts. By incorporating realistic failure scenarios derived from historical tornado data in the Dallas-Fort Worth area, the framework enables proactive, data-informed planning of crew assignments and routing decisions. The integration of multi-agent capacitated vehicle routing and scenario decomposition techniques proved essential in managing computational complexity and ensuring convergence to optimal solutions. Our results demonstrate that the proposed model offers robust and coordinated restoration strategies that can adapt to varying disaster intensities and infrastructure failures.}

\hdk{While the current model assumes fixed crew types and non-overlapping responsibilities, future work may explore more dynamic crew configurations, real-time reassignment strategies, and expanded scenario generation involving critical facility priorities and expert-informed hazard modeling. Additionally, incorporating real-time data streams and operational constraints can further enhance decision support for utility operators during live restoration events. Overall, this study lays a strong foundation for scalable and uncertainty-aware restoration planning, offering both methodological contributions and practical insights for improving power system resilience.}
\newpage


\section*{RESOURCE AVAILABILITY}


\subsection*{Lead contact}


Requests for further information and resources should be directed to and will be fulfilled by the lead contact, Jie Zhang (jiezhang@utdallas.edu).

\subsection*{Materials availability}

This study did not generate new materials.

\subsection*{Data and code availability}
\begin{itemize}
    \item All original code, raw data, and input data have been deposited in Zenodo under \href{https://doi.org/10.5281/zenodo.15786800}{DOI 10.5281/zenodo.15786800} and are publicly accessible. Any additional data referenced in this paper will be provided by the lead contact upon request.
    \item Any further information needed to reanalyze the data presented in this paper is available from the lead contact upon request.
\end{itemize}


\section*{ACKNOWLEDGMENTS}



\hdk{This work was partially supported by the Office of Naval
Research under ONR award number N00014-21-1-2530 and
National Science Foundation under grant 2229417. The United
States Government has a royalty-free license throughout the
world in all copyrightable material contained herein. Any opinions, findings, and conclusions or recommendations expressed in this material are those of the author(s) and do not necessarily reflect the views of the Office of Naval Research and National Science Foundation.}

\section*{AUTHOR CONTRIBUTIONS}

Conceptualization, H.D.K., R.A.J., S.C., and J.Z.; methodology,  H.D.K. and R.A.J.; investigation, H.D.K.  and R.A.J.; writing--original draft, H.D.K.; writing--review \& editing, H.D.K., R.A.J., S.C., and J.Z.; funding acquisition, J.Z.; resources,  H.D.K. and R.A.J.; supervision, S.C. and J.Z.


\section*{DECLARATION OF INTERESTS}


The authors declare no competing interests.


\bibliography{Reference.bib}
\end{document}